\documentclass{amsproc}

\newtheorem{theorem}{Theorem}[section]
\newtheorem{lemma}[theorem]{Lemma}

\theoremstyle{definition}

\theoremstyle{remark}
\newtheorem{remark}[theorem]{Remark}

\numberwithin{equation}{section}

\usepackage{amssymb}
\usepackage{graphicx}
\input xy
\xyoption{all}

\def\R{\mathbb{R}}
 \def\Z{\mathbb Z}
 \def\C{\mathbb C}

\def\sI{\mathcal I}
\def\sT{\mathcal T}

\def\sL{\mathcal L}

\def\MCG{\mbox{MCG}}

\def\P{\mathbb{P}}

 \def\qed{\hfill\framebox(5,5){}}
\DeclareMathAlphabet{\mathpzc}{OT1}{pzc}{m}{it}

\renewcommand{\bold}[1]{\smallskip \noindent {\bf \boldmath #1 }\nopagebreak[4]}

\title{Generalized Lantern Relations and Planar Line Arrangements}
\author{Eriko Hironaka}
\address{Department of Mathematics, Florida State University, Tallahassee, FL 32306-4510}
\email{hironaka@math.fsu.edu}
\thanks{This work was partially supported by a grant from the Simons Foundation (\#209171to Eriko Hironaka).}
\subjclass{Primary 57M27, 20F36;  Secondary 32Q55}
\date{July 1, 2011}
\begin{document}

\begin{abstract}  
In this paper we show that 
to each planar line arrangement defined over the real numbers, for which no two lines are parallel,
one can write down a corresponding
relation on Dehn twists that can be read off from the combinatorics and relative locations of 
intersections.   This gives an alternate proof of Wajnryb's generalized lantern relations,
and of Endo, Mark and Horn-Morris'  daisy relations.
\end{abstract}
\maketitle

\section{Introduction}

Braid monodromy is a useful tool for studying the topology of complements
of line arrangements as is seen in work of \cite{MT:Lines}, 
\cite{Cohen-Suciu97}, \cite{Cordovil98}.
In this paper, we adapt braid monodromy techniques to generate
relations on Dehn twists in the mapping class group $\MCG(S)$ of an oriented
surface $S$ of finite type.

The study of  hyperplane arrangements  has a
rich history in the realms of topology, algebraic geometry, and analysis
(see, for example, \cite{O-T:Arr} for a survey).  
While easy to draw, the deformation theory  
of real planar line arrangements $\sL$ holds many mysteries.  For example, 
there are topologically distinct real line arrangements with equivalent combinatorics
 \cite{Artal04} \cite{Ryb:Fund}.   Moreover, by
 the Silvester-Gallai theorem \cite{Sylvester} there are planar line arrangements
defined over complex numbers, whose combinatorics cannot be duplicated by
a real line arrangement, for example, the lines joining flexes of a smooth cubic
curve.  Braid monodromy is a convenient tool for encoding
the local and global topology of $\sL$.

The lantern relation on Dehn  twists is of special interest because it
with four other simple to state relations generate all relations
in the Dehn-Lickorish-Humphries presentation of $\MCG(S)$ (see \cite{Birman:SantaCruz}, \cite{FM:MCG}, \cite{Matsumoto02},  \cite{Wajnryb:MCG}).
The lantern relation also plays an important role in J. Harer's proof
that the abelianization of $\MCG(S)$ is trivial if $S$ is a closed surface of genus
$g \ge 3$ \cite{Harer83} (cf., \cite{FM:MCG}, Sec. 5.1.2).

\begin{figure}[tb]
\includegraphics[width=3in]{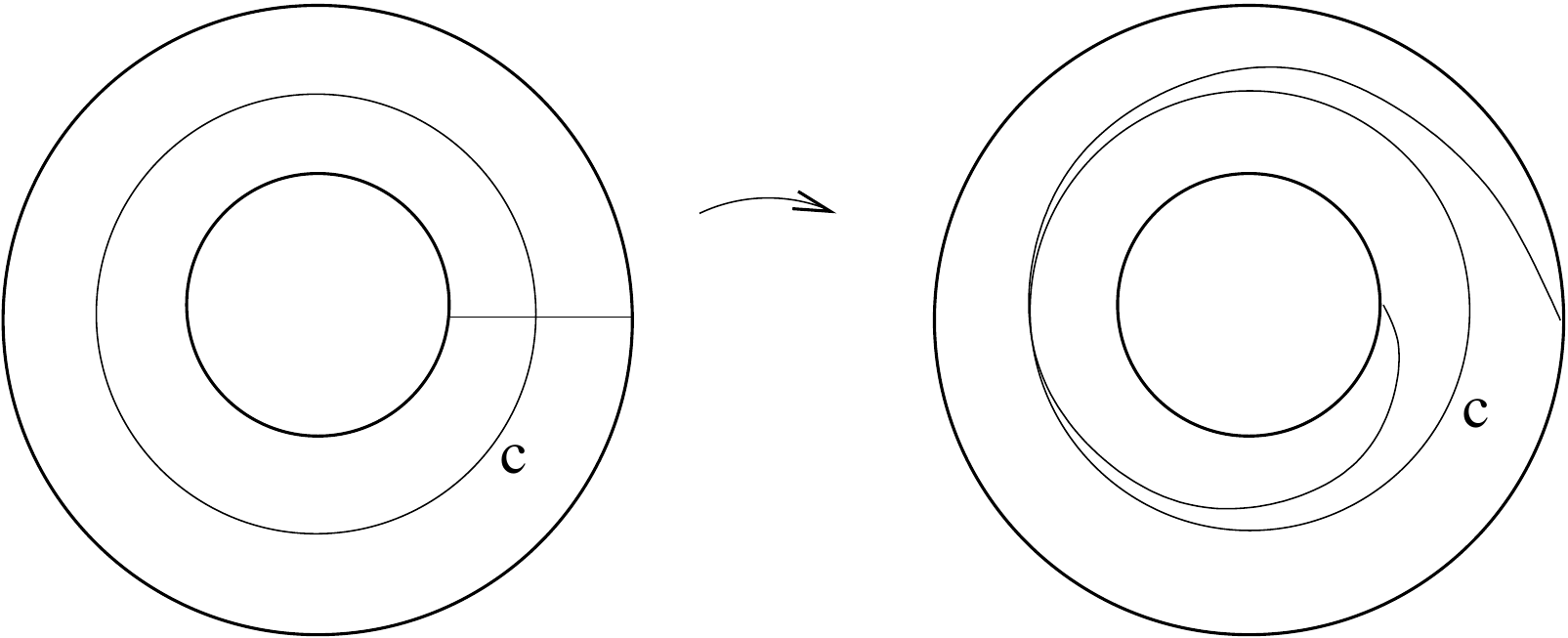}
\caption{Right Dehn twist on an annulus $A$.}
\label{Dehntwist-fig}
\end{figure}

Let $S$ be an oriented surface
of finite type.  If $S$ is closed, the mapping class group $\MCG(S)$ is the group of isotopy classes of
self-homeomorphisms of $S$.  
If $S$ has boundary components,  then the definition of $\MCG(S)$ has the additional
condition that all maps fix the boundary of $S$
pointwise.  For  a compact annulus $A$,  $\MCG(A)$ is isomorphic to $\Z$ and is generated by a right or left Dehn twist
centered at its core curve.  As illustrated in Figure~\ref{Dehntwist-fig}, a right Dehn twist  takes
an arc on $A$ transverse to the core curve to an arc  that wraps
once around the core curve turning in the right hand direction (a left Dehn twist correspondingly turns in the left direction)
as it passes through $c$. 
A Dehn twist can also be thought of as rotating one of the boundary components
by $360{}^\circ$ while leaving the other boundary component fixed.
Each simple closed curve $c$ on $S$ determines a right Dehn twist on an annulus neighborhood of $c$, 
and this Dehn twist extends by the identity to all of $S$.  The isotopy class $\partial_c$ of
this map is the {\it (right) Dehn twist centered at $c$} and is an element of $\MCG(S)$.

The original statement and proof of the lantern relation appears in Dehn's  1938 paper 
 \cite{Dehn38} and relates a product of three interior Dehn twists to four boundary twists
 on a genus zero surface with four boundary components.  
 The relation was  rediscovered by
D. Johnson  \cite{Johnson79}, and  B. Wajnryb gave the following generalized version in 
 \cite{Wajnryb98} (Lemma 17).

\begin{theorem}[Wajnryb]\label{Wajnryb-thm}Let $S^c_{0,n+1}\subset S$ be a surface
of genus zero with $n+1$ boundary components $d_0,d_1,\dots,d_n$.
There is a collection of simple closed curves 
$a_{i,j}, 1\leq i < j \leq n$  in the interior of $S^c_{0,n+1}$, so that
\begin{enumerate}
\item for each $i,j$, $a_{i,j}$ separates
$d_i \cup d_j$ from the rest of the boundary components, and
\item there is a relation on Dehn twists
\begin{align}\label{lantern-eqn}
\partial_0 &(\partial_1\cdots\partial_n)^{n-2}
=\alpha_{1,2} \cdots \alpha_{1,n}\alpha_{2,3}\cdots\alpha_{2,n}\cdots
\alpha_{n-2,n-1} \alpha_{n-2,n}\alpha_{n-1,n},
\end{align}
where $\alpha_{i,j}$ is the right Dehn twist centered at $a_{i,j}$,
and $\partial_i$ is the right Dehn centered at a curve parallel to the boundary components
$d_i$.    
\end{enumerate}
\end{theorem}

We now generalize Theorem~\ref{Wajnryb-thm} in terms of  line arrangements in $\R^2$.

\begin{theorem} \label{Main-thm} Let $\sL$ be a union  of
$n \ge 3$ distinct lines in the $(x,y)$-plane over the reals with distinct slopes and
no  slope parallel to the $y$-axis.
Let $\sI = \{p_1,\dots,p_s\}$ be the intersection points on $\sL$ numbered by largest to smallest
$x$-coordinate.   For each $L \in \sL$,  let $\mu_L$ be the number of points in $\sI \cap L$.
Let $S^c_{0,n+1}$ be a surface of genus zero and $n+1$ boundary components, 
one denoted $d_L$ for each $L \in \sL$, and one extra boundary component $d_0$.
Then there are 
simple closed curve $a_{p_k}$, $k=1,\dots,s$ on $S^c_{0,n+1}$ so that the following holds:
\begin{enumerate}
\item each $a_{p_k}$ separates $$
\bigcup_{p_k \in L \in \sL} d_L 
$$
  from the rest of the boundary curves; and
\item
 the Dehn twist $\partial_L $ centered at $d_L$ and $\alpha_{p_k}$
centered at $a_{p_k}$
satisfy
\begin{eqnarray}\label{lantern1-eqn}
\partial_0 \prod_{L \in \sL} \partial_L^{\mu_L-1} & =& \alpha_{p_s} \cdots \alpha_{p_1}.
\end{eqnarray}
\end{enumerate}
 \end{theorem}

\begin{remark} In Equation~$(\ref{lantern1-eqn})$, the terms on the left side commute,
while the ones on the right typically don't.  Thus, the ordering of $p_1,\dots,p_s$
matters, and reflects the global (as opposed to local) combinatorics of the line arrangement. 
The curves $a_{p_k}$ can be found explicitly (see Section~\ref{twists-sec}, Lemma~\ref{monodromy-lem}).
\end{remark}

\begin{remark} The relations in $\MCG(S^c_{0,n+1})$ give rise to relations on $\MCG(S)$ for any 
surface $S$ admitting an embedding $S^c_{0,n+1} \hookrightarrow S$.
\end{remark}

\begin{figure}[tb]
\includegraphics[width=5in]{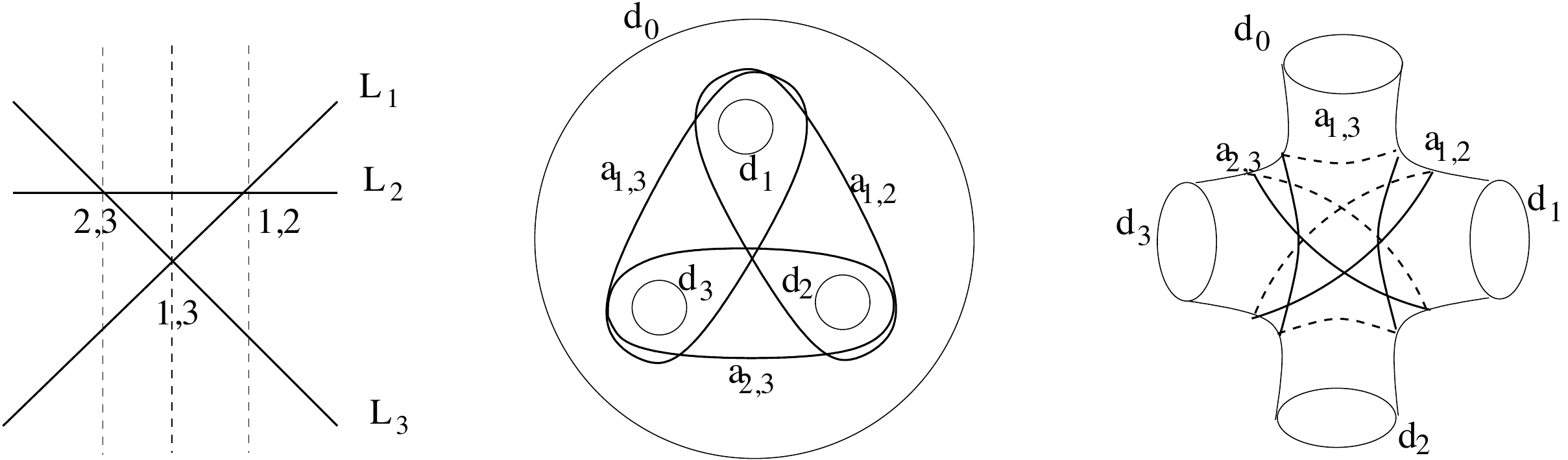}
\caption{Three lines in general position, and curves defining associated lantern relation drawn two ways.}
\label{lantern-fig}
\end{figure}

When $n=3$, Theorem~\ref{Wajnryb-thm} gives the standard Lantern relation
$$
\partial_0\partial_1\partial_2\partial_3 = \alpha_{1,2}\alpha_{1,3}\alpha_{2,3}.
$$
The core curves for these Dehn twists and the corresponding line arrangements
are shown in Figure~\ref{lantern-fig}.  The diagram to the right is the motivation for
the name of this relation.

Here is a sketch of our proof of
 Theorem~\ref{Main-thm}.  First consider a great
ball $B \subset \C^2$ containing all the points of intersection of $\sL$.
Let $\C\P^2$ be the projective compactification of $\C^2$.   Then the
complement of $B$ in $\C\P^2$ is a neighborhood of the ``line at infinity"
or $L_\infty = \C\P^2 \setminus \C^2$.   Let $\rho : \C^2 \rightarrow \C$
be a generic projection.
The monodromy of $\rho$ over the boundary $\gamma$ of a large disk in $\C$
 depends  only on the way  $\sL$ intersects $L_\infty$.
If no lines in $\sL$ are parallel to each other, then 
it is possible to move the lines in $\sL$  to obtain a new configuration
$\sT$ where all lines meet at a single
point  without changing any
slopes, and hence the topology of $\C\P^2 \setminus B$ remains the same (Lemma~\ref{deformation-lem}).     Thus, the monodromies over $\gamma$ defined by
$\sL$ and $\sT$ are the same.     Theorem~\ref{Main-thm} then follows from a description of
 the monodromy of line arrangements on compactified fibers of
a generic projection (Lemma~\ref{monodromy-lem}).   The monodromy can be interpreted as point pushing
maps, where we keep track of twisting on the boundary components of the 
compactified fibers using the complex coordinate system of the ambient
space $\C^2$ (Lemma~\ref{twist-lem}).

This paper is organized as follows.
  In Section~\ref{braidmonodromy-sec} we recall the Moishezon-Teicher braid monodromy representation of a free group associated to a planar line arrangement.    We refine the   representation so that its image is the
  the mapping class group of compactified fibers in Section~\ref{twists-sec}.
  In Section~\ref{def-sec}, we prove Theorem~\ref{Main-thm} using deformations
  of line arrangements and give further variations of
  the lantern relation, including the daisy relation (Theorem~\ref{daisy-thm}).
  
{\bf Acknowledgments:}   The author is grateful to J. Mortada and D. Margalit for helpful discussions and comments, and to the referee for careful corrections to the original version.

  \section{Real line arrangements and relations on Dehn twists}

In this section, we analyze line arrangements $\sL$ in the complex plane defined by
real equations and the monodromy on generic fibers under linear projections
$$
\C^2 \setminus \sL \rightarrow \C.
$$
A key ingredient  is  B. Moishezon and M. Teicher description of
the monodromy  as elements of the pure braid group. (See, for example, \cite{MT:Lines} and \cite{Hironaka:Mem}.)
We generalize this braid monodromy by studying the action of the monodromy  
not only on generic fibers of $\rho$, but also
on their compactifications as genus zero surfaces with boundary.  This leads to a proof
of Theorem~\ref{Main-thm}.

The ideas in this section can be generalized to more arbitrary plane curves.  An investigation 
of the topology of plane curve complements using such general projections
appears in work of O. Zariski and E. van Kampen  \cite{Kampen:Fund}.   We leave this
as a topic for future study.

\subsection{Braid monodromy defined by planar line arrangements over the reals}\label{braidmonodromy-sec}

In this section we recall the Moishezon-Teicher braid monodromy associated to a real line arrangement.
For convenience we choose Euclidean coordinates
$(x,y)$ for $\C^2$ so that no line is parallel to the $y$-axis, and no two intersection
points have the same $x$-coordinate.  For $i=1,\dots,n$, let $L_i$ be the zero set
of a linear equation in $x$ and $y$ with real coefficients:
$$
 L_i = \{(x,y) \ ; \  y=m_i x + c_i\} \qquad m_i,c_i \in \R
$$
and assume that the lines are ordered so that the slopes satisfy:
$$
m_1 > m_2 > \cdots > m_n.
$$
Let $\sI = \sI (\sL) = \{p_1,\dots,p_s\} \subset \C^2$ be the collection of intersections points of  $\sL$
ordered so that the $x$-coordinates are
strictly decreasing.

Let $\rho : \C^2 \rightarrow \C$ be the projection of $\C^2$ onto $\C$ given by $\rho(x,y) = x$.
For each $x \in \C$, let
$$
F_x = \rho^{-1}(x) \setminus \sL.
$$
The $y$-coordinate allows us to uniformly identify $F_x$ with the complement in 
$\C$ of $n$ points $L_i(x)$, where
$$
\{(x,L_i(x))\} =  \rho^{-1}(x) \cap L_i.
$$
Thus, we will think of $F_x$ as a continuous family of copies of $\C$ minus a finite set of points, 
rather than as a subset of $\C^2$.

Let $x_0 \in \R$  be greater than any point in  $\rho(\sI)$.  
Then there is a natural map 
$$
\gamma: [0,1] \rightarrow \C \setminus \rho(\sI)
$$
from arcs 
based at $x_0$ to a braid on $n$ strands in $\C$ parameterized by
$$
\{L_i(\gamma(t)) \ : \ i=1,\dots,n \}.
$$
Since two homotopic
arcs give rise to isotopic braids, and a composition of arcs gives rise to a composition of braids, we have  a homomorphism
$$
\beta : \pi_1(\C \setminus \rho(\sI),x_0) \rightarrow B(S^2,n+1)
$$
from the fundamental group to the spherical braid group on $n+1$ strands.

 The {\it (braid) monodromy} of $(\C^2, \sL)$ with respect to the projection $\rho$
and basepoint $x_0$ is the homomorphism
\begin{eqnarray}
\sigma_{\sL} : \pi_1(\C \setminus \rho(\sI),x_0) \rightarrow \MCG(F_{x_0}),
\end{eqnarray}
given by the composition of $\beta$
and the  braid representation
$$
B(S^2,n+1) \rightarrow \MCG(S_{0,n+1}) = \MCG(F_{x_0}),
$$
from the braid group to the mapping class group on a genus zero surface with $n+1$ punctures.

We now study the image of simple generators of $\pi_1(\C \setminus \rho(\sI),x_0)$
 in $\MCG(F_{x_0})$.
\begin{figure}[tb]
\includegraphics[width=3.5in]{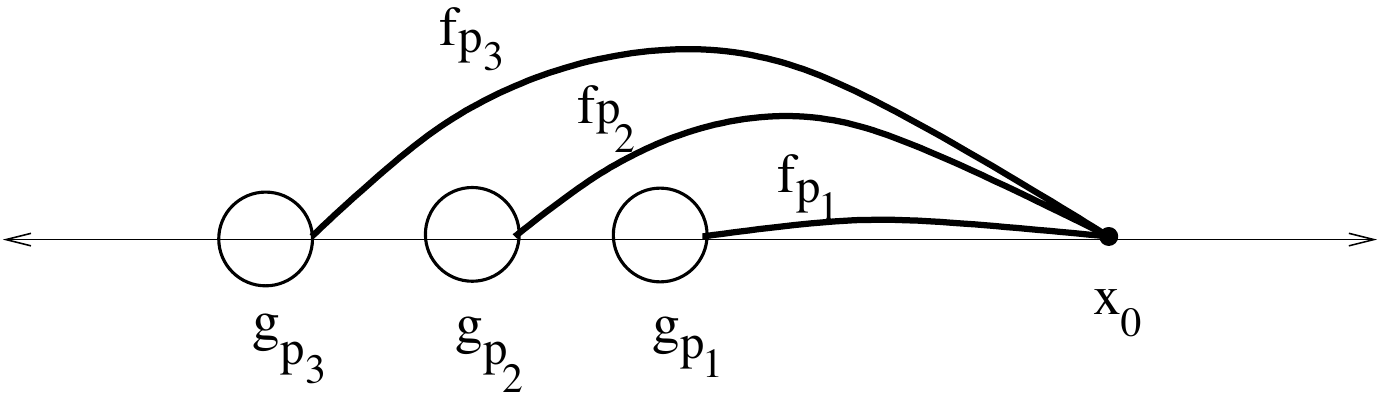}
\caption{Simple loop generators for $\pi_1(\C \setminus \rho(\sI))$.}
\label{firstfig}
\end{figure}

By a simple loop in $\pi_1(\C \setminus \rho(\sI),x_0)$,
we mean a arc of the form $\ell_p = f_p g_p f_p^{-1}$, where $p \in \rho(\sI)$, $\epsilon_p > 0$,
\begin{eqnarray*}
g_p : [0,1] &\rightarrow& \C \setminus \rho(\sI)\\
t &\mapsto& p + \epsilon_p e^{2\pi i t}
\end{eqnarray*}
and $f_p$ is a arc from $x_0$ to $p+{\epsilon_p}$ whose image is in the upper half plane
except at its endpoints.  Since $\pi_1(\C \setminus \rho(\sI),x_0)$ is generated by
simple loops, it is enough to understand the monodromy in the image of these elements.

In order to describe the monodromy of $\ell_p$ we 
study how $F_x$ is transformed as $x$ follows its arc segments $g_p$ and $f_p$
First we look at $g_p$.  Let $L_{j_1},L_{j_2},\dots,L_{j_k}$ be the lines in $\sL$ that pass through
$p$.   
We can assume by a translation of coordinates that $p = 0$, and $L_{j_r}$ is defined by an equation of the
form
$$
y = m_{r} x
$$
where $m_1 > m_2 > \cdots > m_k$.   Then as $t$ varies in $[0,1]$, the intersection of $L_{j_r}$ with $\rho^{-1}{(g_p(t))}$
is given by
$$
L_{j_r}(g_p(t)) = (\epsilon_p e^{2\pi i  t}, m_r \epsilon_p e^{2\pi i t}).
$$

\begin{figure}[tb]
\includegraphics[width=3in]{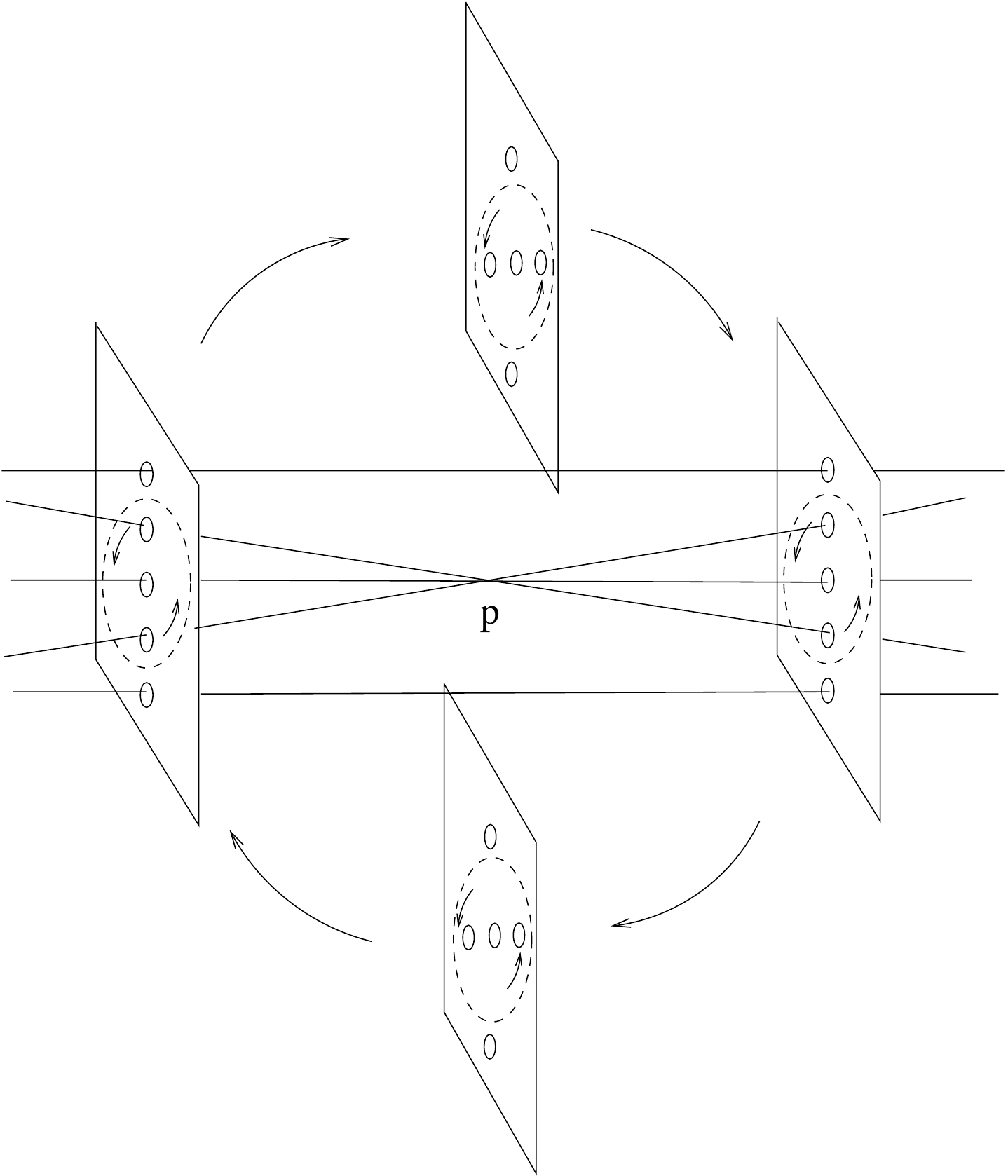}
\caption{Monodromy defined by $g_p$ with the real part of $\sL$ drawn in.}
\label{monodromy-fig}
\end{figure}

The other lines in $\sL$ locally can be thought of as having constant slope, hence their intersections with $\rho^{-1}({g_p(t)})$
retain their order and stay outside a circle on $F_{g_p(t)}$
enclosing $L_{j_1}(g_p(t)),\dots,L_{j_k}(g_p(t))$  (see Figure~\ref{monodromy-fig}).  Let $a_p^{\mbox{\small loc}} \subset F_{p + \epsilon}$ be this circle.
The restriction of $\rho$ to $\C^2 \setminus \sL$ defines a trivial bundle over the image of $f_p$.  Thus
  $a_p^{\mbox{\small loc}}$ determines a simple closed curve $a_p$ on $F_{x_0}$ separating 
 $L_{j_1}(x_0),\dots,L_{j_k}(x_0)$ from the rest of the $L_j(x_0)$.

Next we notice that lifting over $f_p$ defines a mapping class on $F_{x_0}$.    This is because 
there is a canonical identification of $F_{x_0}$ and $F_x$ for
any $x \in \R \setminus \rho(\sI)$ given by the natural ordering of $\sL \cap \rho^{-1}(x)$ 
by the size of the $y$-coordinate from largest to smallest.
Thus $f_p$ determines a braid on $n$ strands and corresponding mapping class
$\beta_p \in \MCG(F_{x_0}^c)$.    We have shown the following.

\begin{lemma}\label{monodromy-lem} Let $\ell_p =  f_pg_pf_p^{-1}$.  The element $\sigma_\sL(\ell_p)$ in  $\MCG(F_{x_0}^c)$
is the Dehn twist $\alpha_p$ centered at $a_p = \beta_p^{-1}(a_p^{\mbox{loc}})$.
\end{lemma}

{\bf Proof.}   By the above descriptions of the fibers above the arcs $f_p$ and $g_p$, we
 can decompose $\sigma_{\sL}(\ell_p)$ as
$$
\alpha_p = \beta_p^{-1} \circ \sigma_p \circ \beta_p,
$$
where $\sigma_p$ is a right Dehn twist centered at $a_p^{\mbox{loc}} = \beta_p(a_p)$.  \qed

\subsection{Monodromy on compactified fibers}\label{twists-sec}

In this section, we define the monodromy representation of $\pi_1(\C^2 \setminus \sL,y_0)$
into $\MCG(F_{x_0}^c)$, where $F_{x_0}^c$ is a compactification of $F_{x_0}$ as a
compact surface with boundary.  

As before choose coordinates for $\C^2$, and 
let $\sL = \cup_{i=1}^n L_i$ be a planar line arrangement defined over the reals with
distinct slopes.   Assume all points of intersection $\sI$ have distinct $x$-coordinates.
Let $\epsilon > 0$ be such that the $\epsilon$ radius disks $N_\epsilon(p)$ around the points 
$p \in \rho(\sI)$ are disjoint.   Let $\delta > 0$ be such that the $\delta$ radius tubular neighborhoods
$N_\delta(L_i)$
 around $L_i$ are disjoint in the complement of 
$$
 \bigcup_{p \in \rho(\sI)} \rho^{-1}(N_\epsilon(p)).
$$
Let $D$ be a disk in $\C$ containing all points of $\rho(\sI)$ in its interior, and having
$x_0$ on its boundary.
Let $N_\infty$ be  a disk centered at the origin of $\C$ so that $\C \times N_\infty$
contains $\sL \cap \rho^{-1}(D)$.

For each $x \in \C \setminus \rho(\sI)$, let 
$$
F_x^c = \rho^{-1}(x) \cap \left (\C \times N_\infty \setminus \cup_{i=1}^n N_\delta(L_i)\right ) \quad \subset \quad F_x \setminus N_\delta(L_i).
$$
For each $x \in D$ and $i=1,\dots,n$, let 
$$
d_i(x) = \partial N_\delta(L_i) \cap F_x.
$$

We are now ready to define the monodromy on the compactified fibers
$$
\sigma^c_{\sL}: \pi_1(F_{x_0},y_0) \rightarrow \MCG(F^c_{x_0}).
$$
Let $\eta$ be the  {\it inclusion homomorphism}
$$
\eta : \MCG(F_{x_0}^c)  \rightarrow \MCG(F_{x_0}),
$$
that is, the homomorphism induced by inclusion $F_{x_0}^c \subset F_{x_0}$.
Then we would like to have a commutative diagram
$$
\xymatrix{
\pi_1(\C \setminus \rho(\sI),y_0)\ar[r]^{\sigma^c_{\sL}}\ar [dr]^{\sigma_{\sL}}& \MCG(F_0^c)\ar[d]^\eta\\
& \MCG(F_0).
}
$$
The kernel of $\eta$ is generated by 
Dehn twists centered at the boundary components of $F^c_{x_0}$
(Theorem 3.18, \cite{FM:MCG}).
Thus, in order to describe $\sigma^c_{\sL}$, we need to understand what twists occur near boundary
components in the monodromy associated to the arcs $g_p$
and $f_p$ defined in Section~\ref{braidmonodromy-sec}.

Consider the simplest case when  $L \subset \C^2$ is
a single line defined by
$$
y= L(x) = mx.
$$
Let $N_\delta(L)$ be the tubular neighborhood around $L$ in $\C^2$
$$
N_\delta(L) = \{(x,L(x)+y) \ : \ |y| < \delta\}.
$$
Then $N_\delta(L) \cap F_{g(t)}$ is a disk centered at $L(g(t))$ of
radius $\delta$.   The boundary $\partial N_\delta(L)$ is a trivial
bundle over $\C \setminus  \delta(\rho(\sI))$ with trivialization defined
by the framing of $\C$ by real and purely imaginary coordinates.

Now assume that there are several lines $L_{j_1},\dots,L_{j_k}$ meeting above $p \in \rho(\sI)$.
Let
$L$ be a line through $p$ with slope equal to the average of those of $L_{j_1},\dots,L_{j_k}$,
and let $\delta >0$ be such that $N_\delta(L) \cap F_{g_p(t)}$ contains
$d_{j_1}(g_p(t)), \dots, d_{j_k}(g_p(t))$, but no other boundary components
of $F^c_{g_p(t)}$, for all $t$.   Let $d_p(t)$ be the boundary component of $F_{g_p(0)}$ given by
$$
d_p(t) = \partial N_\epsilon(L) \cap F_{g_p(t)}.
$$
Let $d_p = d_p(0)$ and $d_{j_i} = d_{j_i}(g_p(0))$.

Then looking at
Figure~\ref{monodromy-fig}, we see that the points $L_{j_1}(g_p(t)),\dots, L_{j_k}(g_p(t))$ rotate as a group $360{}^\circ$ in
the counterclockwise direction as $t$ ranges in $[0,1]$.    The corresponding mapping class on the bounded portion of
$F_{g_p(0)}$ enclosed by
$d_p$ is the composition of
a clockwise full rotation of $d_p$ and a counterclockwise rotation around $d_{j_1}, \dots, d_{j_k}$.
It can also be thought of as moving the inner boundary components $d_{j_1}(g_p(0))$  in a clockwise direction
while leaving all orientations of boundary components fixed with respect to the complex framing of $\C$.

\begin{figure}[tb]
\includegraphics[width=4in]{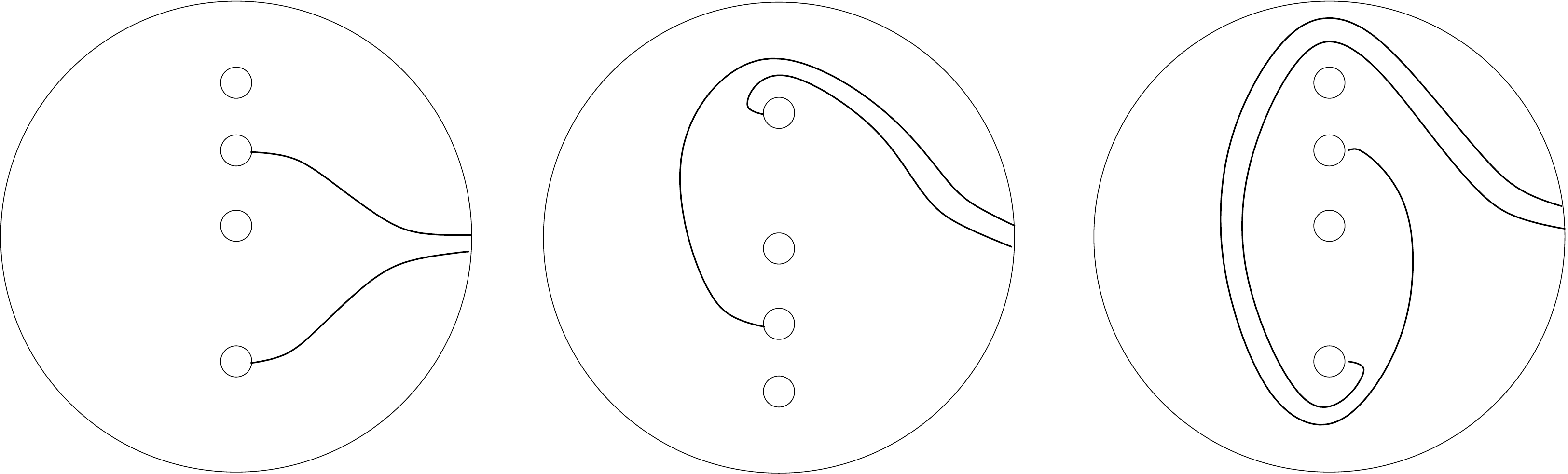}
\caption{The mapping class  $\partial_{d_p}$.}
\label{outerD-fig}
\end{figure}

\begin{figure}[tb]
\includegraphics[width=4in]{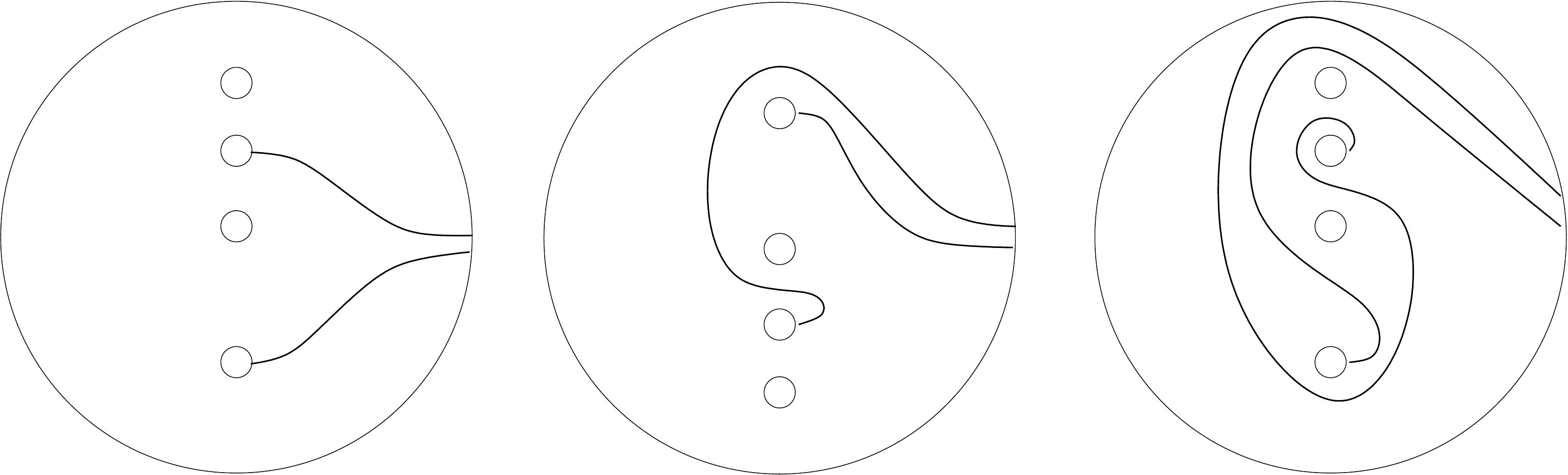}
\caption{The monodromy defined by $g_p$.}
\label{rotation-fig}
\end{figure}

Figure~\ref{outerD-fig} illustrates the Dehn twist $\partial_{d_p}$ centered 
 at a simple closed curve parallel to $d_p$ and Figure~\ref{rotation-fig} shows  the monodromy $\sigma_{\sL}^c(g_p)$ in the case
 when $\sL$ is a union of 4 lines meeting at a single point $p$.  In both figures, the middle picture illustrates the
 fiber $F_{g_p(0.5)}$ half way around the circle traversed by $g_p$.
From the above discussion, we have
$$
\sigma^c_{\sL}(g_p) = (\partial_{d_1}\partial_{d_2}\partial_{d_3}\partial_{d_4})^{-1} \partial_{d_p}.
$$
More generally we have the following lemma.

\begin{lemma}\label{twist-lem}   Let $L_{j_1},\dots,L_{j_k}$ be the lines
meeting above $p$, and let 
\begin{eqnarray*}
g_p : [0,1]& \rightarrow& \C \setminus \sL\\
t &\mapsto& p + \epsilon e^{2\pi i t}.
\end{eqnarray*}
 Then the monodromy on $F_{g_p(0)}^c$ 
defined by $g_p$ is given by  
 $$
 \sigma^c_{\sL}(g_p) = (\partial_{d_{j_1}} \cdots \partial_{d_{j_k}})^{-1} \partial_{d_p}.
 $$
 \end{lemma}

\subsection{Deformations of line arrangements}\label{def-sec}

To finish our proof of Theorem\ref{Main-thm} we analyze the effect of deforming a line arrangement.

Let
$$
\sL=\bigcup_{i=1}^n  L_i
$$
be a finite union of real lines in the Euclidean plane, $\R^2$ with no two lines parallel.
  Let $\sT$ be the
complexified real line arrangement with all $n$ lines intersecting at a single point $p_0$.
Let $\rho : \C^2 \rightarrow \C$ be a generic projection, and let $D \subset \C$ be a disk 
of radius $r$ centered at the origin containing $\rho(\sI)$ and $\rho(p_0)$ in its interior.  Let 
\begin{eqnarray*}
\gamma: [0,1] &\rightarrow& \C \\
t &\mapsto& r e^{2\pi i t}.
\end{eqnarray*}

\begin{lemma}\label{deformation-lem}  
The monodromies $\sigma^c_{\sL}(\gamma)$ and $\sigma^c_{\sT}(\gamma)$
are the same.
\end{lemma}

{\bf Proof.}   Let 
$$
D_\infty = N_\infty \times \C \cap \rho^{-1}(\partial D).
$$
Then $D_\infty \setminus \sT$ and $D_\infty \setminus \sL$ are isomorphic
as fiber bundles over $\gamma$ and hence the monodromies over $\gamma$
defined by $\sL$ and $\sT$ are the same up to isotopy.\qed

\begin{figure}[tb]
\includegraphics[width=5in]{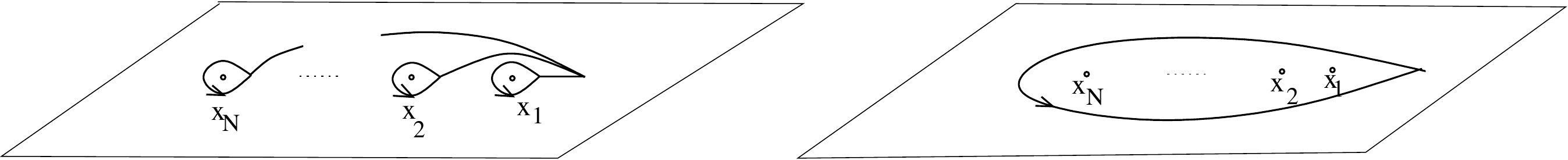}
\caption{Two representatives of $\gamma$ in $\pi_1(\C \setminus \rho(\sI))$.}
\label{pathhomotopy-fig}
\end{figure}

\bold{Proof of Theorem~\ref{Main-thm}.} 
By Lemma~\ref{deformation-lem}, $\sigma_{\sL}^c(\gamma) = \sigma_{\sT}^c(\gamma)$.
Figure~\ref{pathhomotopy-fig} gives an illustration of two equivalent representations
of the homotopy type of $\gamma$.  

By Lemma~\ref{monodromy-lem} and Lemma~\ref{twist-lem}, we have
$$
\sigma^c_{\sT}(\gamma) = (\partial_{d_1} \cdots \partial_{d_n})^{-1}\partial_{d_p}.
$$
Let $p_1,\dots,p_s$ be the elements of $\sI$ numbered by decreasing $x$-coordinate.
Then for each $i=1,\dots,s$, we have 
$$
\sigma_{\sL}^c(f_{p_i}g_{p_i}f_{p_i}^{-1}) = (\partial_{d_{j_1}} \cdots  \partial_{d_{j_k}})^{-1}\alpha_{p_i}
$$
where $\alpha_{p_i}$ 
is the pullback of $d_{p_i}$ along the arc $f_{p_i}$.
Thus, 
$$
\sigma_{\sL}^c(\gamma) = (\partial_{d_1}^{\mu_1}\cdots \partial_{d_s}^{\mu_s})^{-1}\alpha_{p_s}\cdots \alpha_{p_1},
$$
where $\mu_i$ is the number of elements in $\sI \cap L_i$.\qed

 To show that Theorem~\ref{Wajnryb-thm} follows from Theorem~\ref{Main-thm}, 
 we need to show that the ordering given in Equation~$(\ref{lantern-eqn})$
can be obtained by a union of lines $\sL$ satisfying the conditions.
To do this, we start with a union of lines $\sT$ intersecting in a single point.
Let $L_1, \dots L_n$ be the lines in $\sT$ ordered from largest to smallest slope. 
Translate $L_1$ in the positive $x$ direction
 without changing its slope so that 
the intersections of the translated line $L'_1$  with $L_2,\dots,L_n$ have decreasing
$x$-coordinate.  Continue for each line from highest to lowest slope, making sure with each 
time that the shifting $L$ creates new intersections
lying to the left of all previously created ones.      

More generally, we can 
deform the lines through a single point $\sT$ to one in general position $\sL$
so that the only condition on the resulting ordering on the pairs of lines is the
following.  A pair $(i,j)$ must preceed $(i,j+1)$ for each $1 \leq i < j \leq n$.
Thus, we have proved the following restatement of Theorem~\ref{Wajnryb-thm}.

\begin{theorem}\label{order-thm}  Let $\{p_1,\dots,p_s\}$ be an ordering of the pairs $(i,j)$,
$1 \leq i < j \leq n$, so that for all $i$, the sequence
$$
(i,i+1), (i,i+2), \dots, (i,n)
$$
is strictly decreasing.   Then there  a lantern relation of the form
$$
\partial_0 (\partial_1 \cdots \partial_n)^{n-2} = \alpha_{p_1} \cdots \alpha_{p_s}.
$$
\end{theorem}

\section{Applications}

Although it is known that all relations on the Dehn-Lickorish-Humphries generators can
be obtained from the braid, chain, lantern and hyperelliptic relations, there are some other
nice symmetric relations that come out of line arrangements that are not trivially 
derived from the four generating ones.   We conclude this paper with a sampling.

\subsection{Daisy relation}
\begin{figure}[tb]
\includegraphics[width=3.5in]{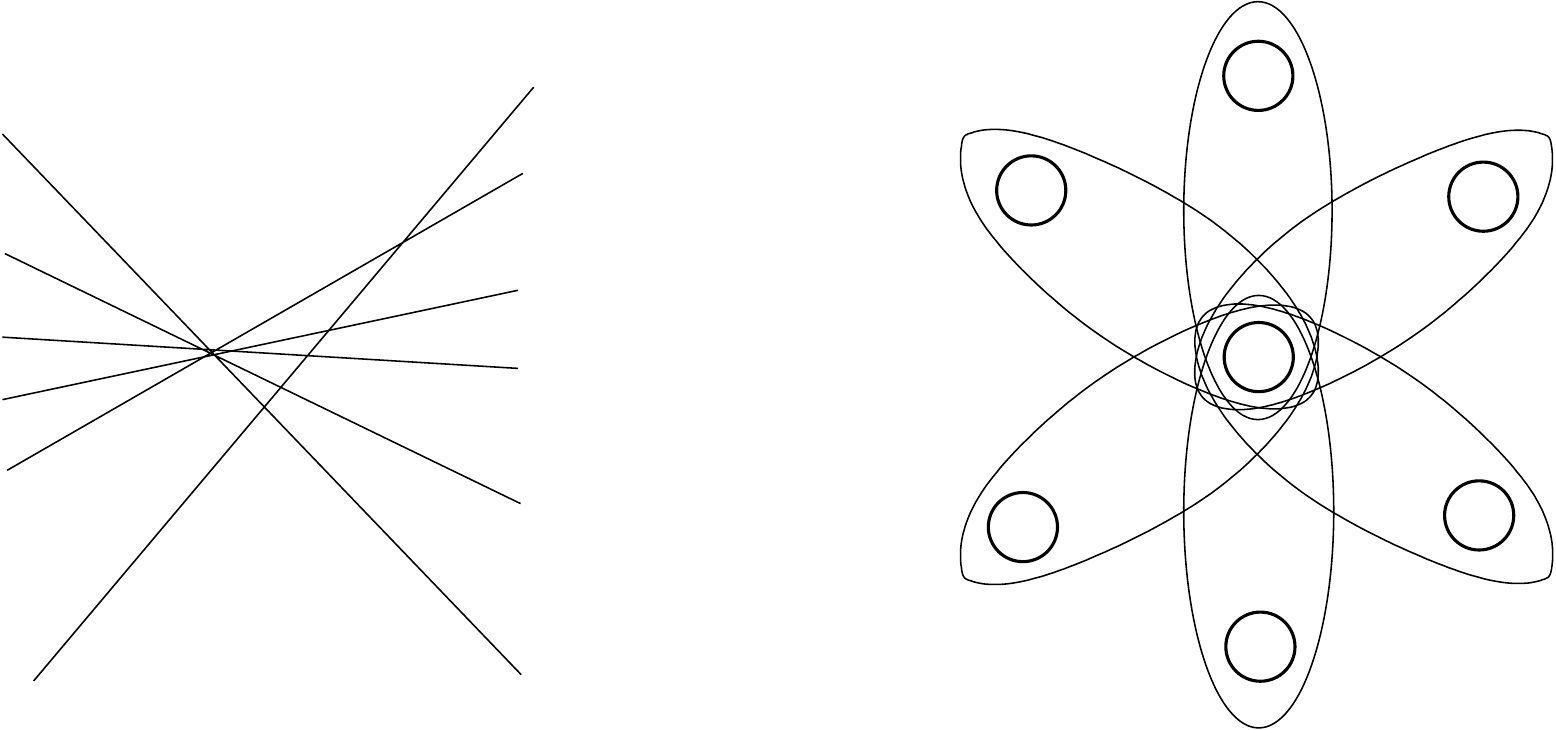}
\caption{Line arrangement, and associated arrangement of curves (n=6).}
\label{starlines-fig}
\end{figure}

Consider the line 
arrangements given in Figure~\ref{starlines-fig}. 
As pointed out to me by D. Margalit, this relation was recently also discovered by
H. Endo, T. Mark, and J. Van Horn-Morris using rational blowdowns of 4-manifolds
\cite{Endo11}.  We follow their nomenclature and call this the  {\it daisy relation}.

Let $S^c_{0,n+1}$ denote the compact surface of genus $0$ with $n+1$ boundary components.
Consider the configuration of simple closed curves shown in Figure~\ref{starlines-fig}.
Let $d_0,\dots,d_n$ be the boundary components of $S^c_{0,n+1}$.
Let $d_1$ be the distinguished boundary component at the center of
the arrangement, and let $d_0, d_2,\dots,d_n$ be the boundary components arranged in 
a circle  (ordered in the clockwise direction around $d_1$).  Let $a_{1,k}$ be 
a simple closed loop encircling $d_1$ and $d_k$, where $k=0, 2, 3,\dots, n$.
Let $\partial_i$ be the Dehn twist centered at $d_i$, and let $\alpha_{1,k}$ be the
Dehn twist centered at $a_{1,k}$.

\begin{theorem}\label{daisy-thm}[Daisy relation]  For $n \ge 3$, the Dehn twists on $S^c_{0,n+1}$
satisfy the relation
$$
\partial_0 \partial_1^{n-2} \partial_2 \cdots \partial_n = \alpha_{1,0}  \alpha_{1,n} \cdots  \alpha_{1,2}
$$
where $\partial_i$ is the Dehn twist centered at  the boundary component $d_i$, and $\alpha_{1,j}$ is the Dehn twist centered at curves
$a_{1,j}$.
\end{theorem}

When $n=3$, Theorem~\ref{daisy-thm} specializes to the usual lantern relation.

{\bf Proof.} We associate the boundary component $d_i$ with $L_i$ for $i=1,\dots,n$,
and $d_0$ with the ``line at infinity".  
Theorem~\ref{Main-thm} applied to the line arrangement in Figure~\ref{starlines-fig}
gives:
$$
\partial_0(\partial_1\cdots \partial_n)^{-1} = R_{p_n} \dots R_{p_1}
$$
where $p_1,\dots,p_n$ are the intersection points of the line arrangement $\sL$
ordered by largest to smallest $x$-coordinate.
For this configuration, $p_k$ gives rise to 
$$
R_{p_k} = (\partial_1\partial_{k+1})^{-1}\alpha_{1,k+1},
$$
for $k=1,\dots,n-1$.  Noting that the loop that separates $d_2 \cup \cdots \cup d_n$
from $d_0 \cup d_1$ can be written as $a_{1,0}$, we have
$$
R_{p_n} = (\partial_2 \cdots \partial_{n-1})^{-1}\alpha_{1,0}
$$
yielding the desired formula.\qed

\begin{figure}[tb]
\includegraphics[width=4in]{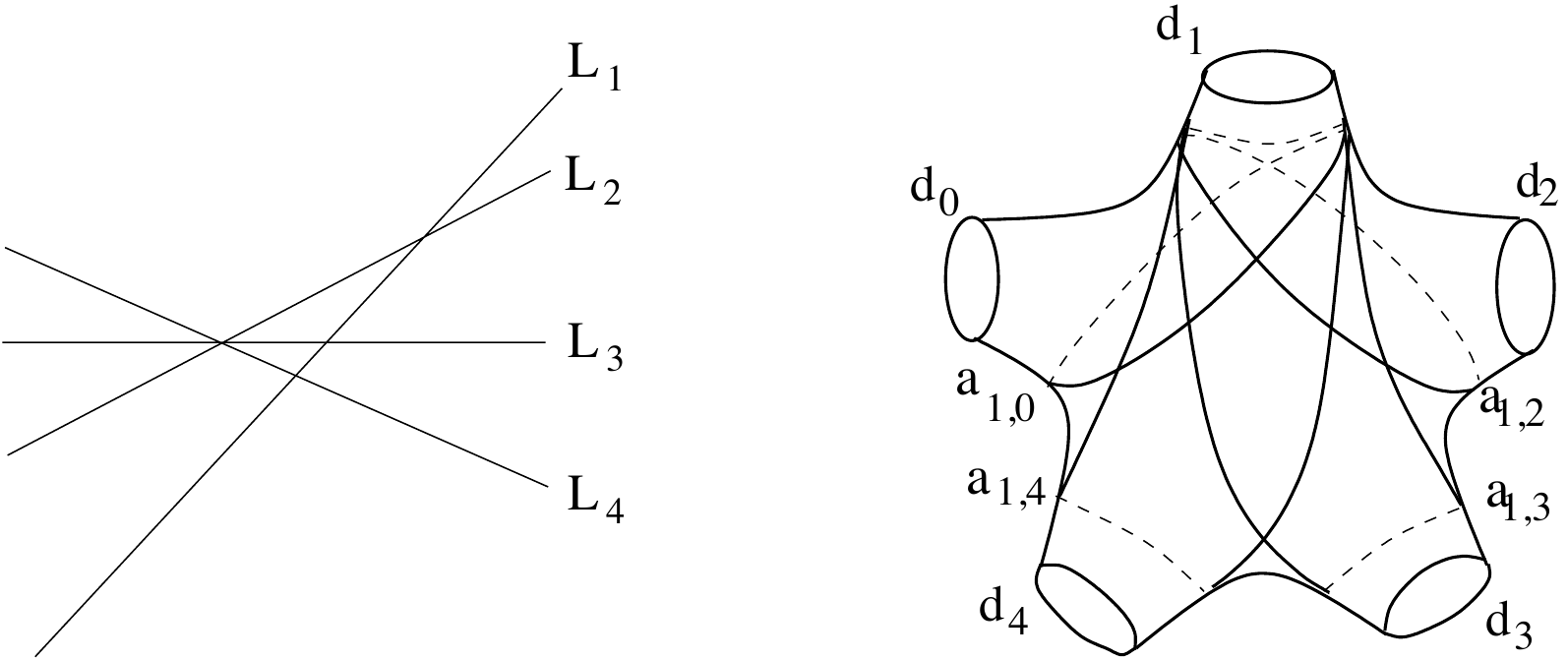}
\caption{Alternate drawing of the daisy configuration (n=4).}
\label{daisy-fig}
\end{figure}

\begin{remark}   Let 
$$
\beta  : B(S^2,n+1) \rightarrow \MCG(S_{0,n+1})
$$
be the braid representation from the spherical braid group to the mapping class group.
Recall the relation $R$ in $B(S^2,n+1)$ given by
\begin{eqnarray*}
(\sigma_1^2) (\sigma_1^{-1}\sigma_2^2\sigma_1)\cdots(\sigma_1^{-1}\sigma_2^{-1}\cdots \sigma_{n-1}^{-1}\sigma_n^2\sigma_{n-1}\cdots\sigma_1)
&=&\sigma_1\cdots \sigma_{n-1} \sigma_n^2 \sigma_{n-1} \cdots \sigma_1.\\
&=&1
\end{eqnarray*}
This induces a relation $R'$ in $\MCG(S_{0,n+1})$.   The daisy relation can be considered
as the lift of $R'$ under the inclusion homomorphism $\eta$.
\end{remark}

\subsection{Doubled daisy relation}

\begin{figure}[tb]
\includegraphics[width=2in]{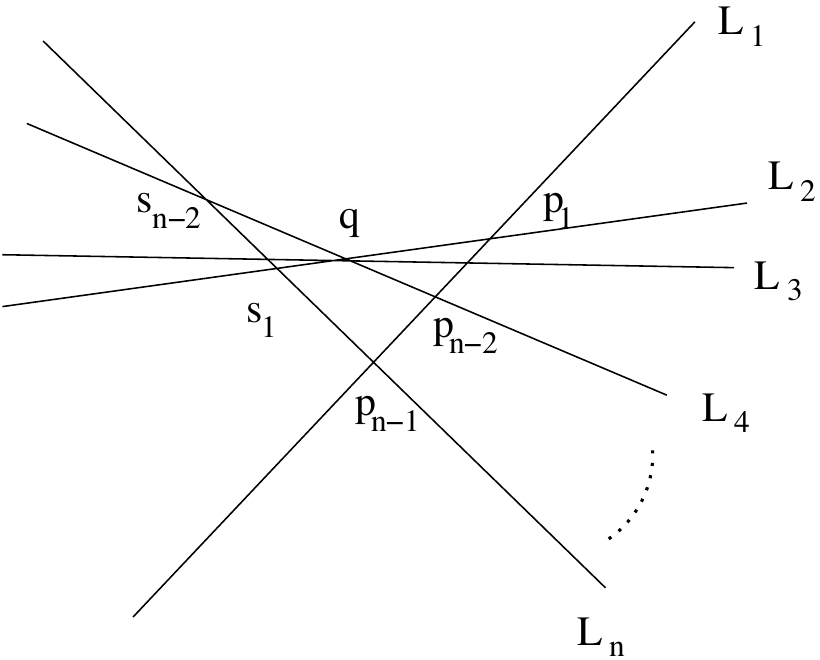}
\caption{Configuration of lines  giving rise to the doubled daisy relation.}
\label{saucerlines-fig}
\end{figure}

As a final example, we consider a configuration
of $n \ge 5$ lines, with $n-2$ meeting in a single point.
There are several ways this can be drawn.  We give one example 
  in Figure~\ref{saucerlines-fig}.  Other line arrangements
  satisfying these conditions will give similar relations, but
  the drawings of the associated curves will be more complicated.

\begin{figure}[tb]
\includegraphics[width=2.25in]{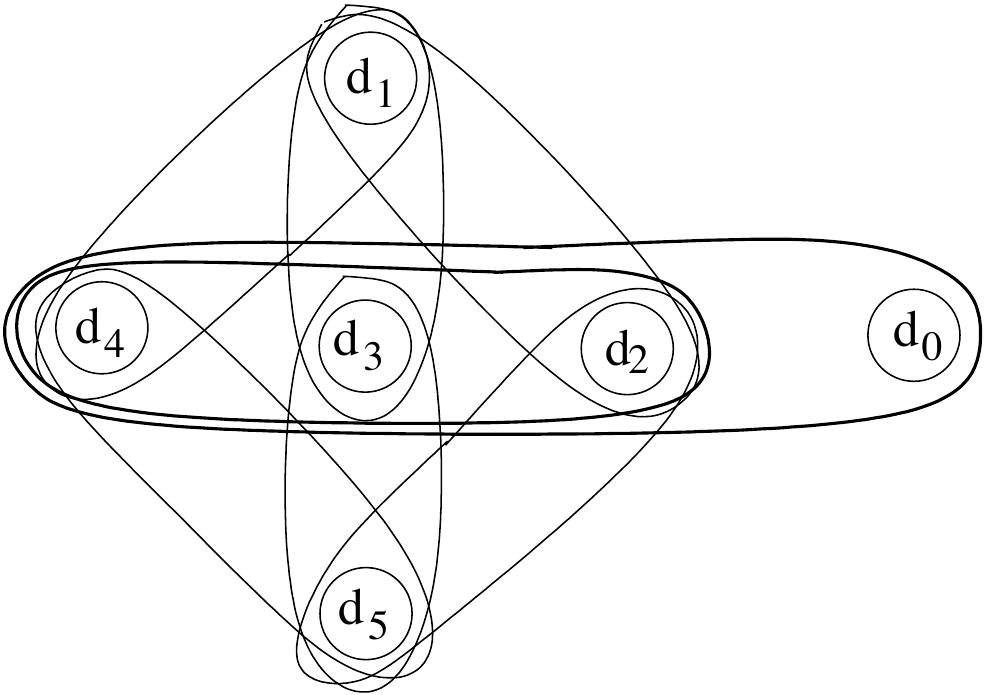}
\caption{The doubled daisy relation for $n=5$.}
\label{saucer5-fig}
\end{figure}

\begin{figure}[tb]
\includegraphics[width=3.75in]{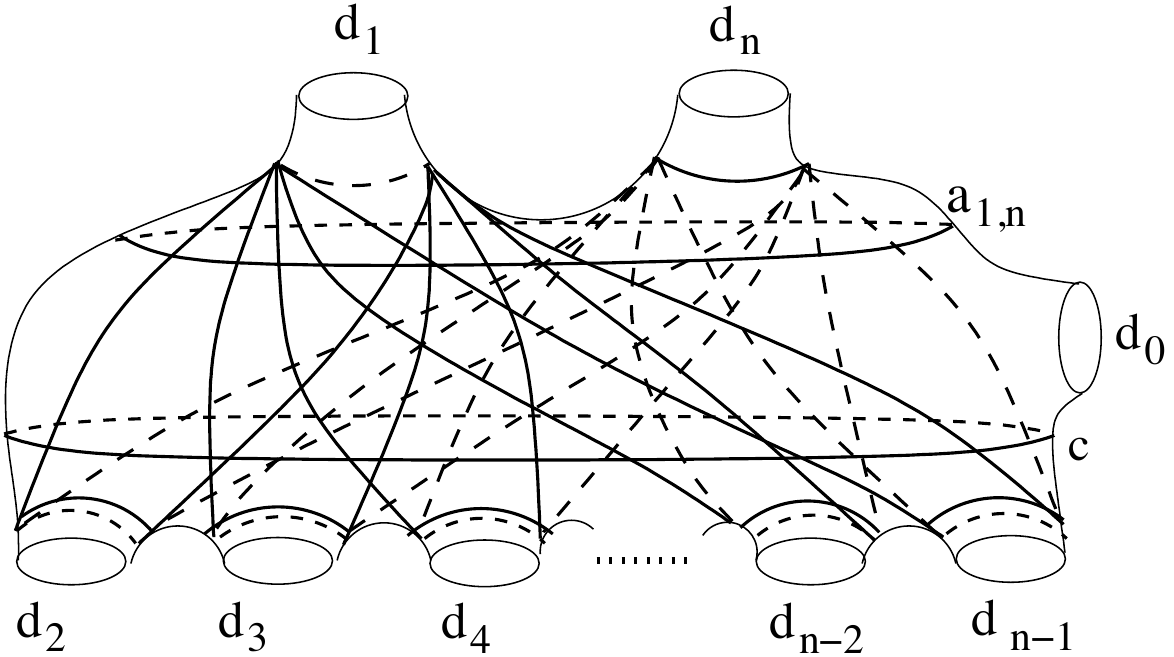}
\caption{Drawing of the general doubled daisy configuration.}
\label{saucer-fig}
\end{figure}

As before, let $d_0,\dots,d_n$ be the boundary components of $S_{0,n+1}^c$.
The boundary component
$d_i$ is associated to the line $L_i$ for $i=1,\dots,n$, and $d_0$ is the boundary
component associated to the ``line at infinity".    Let $a_{i,j}$ be the loop in Figure~\ref{saucer-fig} encircling
$d_i \cup d_j$  and no other boundary component.  Let $c$ be the loop encircling
$d_2,\dots,d_{n-1}$ in Figure~\ref{saucer-fig} (or, when $n=5$,
$d_2, d_3,$ and $d_4$ in Figure~\ref{saucer5-fig}).   

\begin{theorem}[Doubled daisy relation]\label{double-thm} Let $\partial_i$ be the right Dehn twist
centered at $d_i$, $\alpha_{i,j}$ the right Dehn twist centered at $a_{i,j}$, and
$\beta$ the right Dehn twist centered at $c$.
Then
$$
\partial_0  \partial_1^{n-2}  \partial_2 \cdots \partial_{n-1}\partial_n^{n-2}=
\alpha_{n-1,n}\alpha_{n-2,n} \cdots \alpha_{2,n}\ \beta\ \alpha_{1,n} \alpha_{1,n-1} \cdots \alpha_{1,2} 
$$
\end{theorem}

{\bf Proof.}  Theorem 2 applied to  the line arrangement in Figure~\ref{saucerlines-fig} gives
the equation
$$
\partial_0(\partial_1 \cdots \partial_n)^{-1} = R_{s_{n-2}} \cdots R_{s_1} R_q R_{p_{n-2}}
\cdots R_{p_1},
$$
where
\begin{eqnarray*}
R_{p_k} &=& (\partial_1 \partial_{k+1})^{-1}\alpha_{1,k+1}\\
R_q &=& (\partial_2 \cdots \partial_{n-1})^{-1} \beta\\
R_{s_k} &=& (\partial_n \partial_{k+1})^{-1} \alpha_{k+1,n}.
\end{eqnarray*}
(As one sees from Figure~\ref{saucerlines-fig} and Figure~\ref{saucer-fig}, the order of $R_{p_{n-1}}$ and $R_q$
may be interchanged.)
Putting these together yields the desired formula.\qed

\bibliographystyle{amsplain}

\bibliography{math}

 \end{document}